\documentclass{article}

\usepackage{PRIMEarxiv}
\usepackage{amsmath}
\usepackage{subcaption} 

\usepackage[ruled,vlined]{algorithm2e}
\usepackage{float}
\usepackage[utf8]{inputenc} 
\usepackage[T1]{fontenc}    
\usepackage{hyperref}       
\usepackage{url}            
\usepackage{booktabs}       
\usepackage{amsfonts}       
\usepackage{nicefrac}       
\usepackage{microtype}      
\usepackage{lipsum}
\usepackage{fancyhdr}       
\usepackage{graphicx}       
\graphicspath{{media/}}     

\title{Noise Removal in One-Dimensional Signals using Iterative Shrinkage Total Variation Algorithm
}

\author{
  Joyce Oliveira dos Santos \\
  Federal University of Rio Grande do Norte \\
  Caicó- RN\\
\texttt{joyce.santos.709@edu.ufrn.br} \\
   \And
  Francisco Márcio Barboza \\
  Department of Computing and Technology  \\
  Federal University of Rio Grande do Norte\\
  Caicó- RN\\
  \texttt{marcio.barboza@ufrn.br} \\
}

\begin{document}
\maketitle

\begin{abstract}
The total variation filtering technique emerges as a highly effective strategy for restoring signals with discontinuities in various parts of their structure. This study presents and implements a one-dimensional signal filtering algorithm based on total variation. The aim is to demonstrate the effectiveness of this algorithm through a series of synthetic filtering tests. The results presented in this paper were significant in demonstrating the proposed algorithm's effectiveness. Through a series of rigorously conducted experiments, the algorithm's ability to solve complex noise removal problems in various scenarios was evidenced.
       
\end{abstract}


\section{Introduction}
The Total Variation (TV) Filtering technique was designed to reduce noise while preserving sharp edges in the original signal. This approach was introduced by \cite{rudin1992nonlinear}. Noise removal by the TV 'filter' is achieved by minimizing a specific cost function \cite{selesnick2012total}.

Various algorithms have been proposed to implement TV filtering, including the one described by \cite{selesnick2010total}. Although the algorithm can be derived in several ways, the derivation presented here is based on descriptions from specific sources. TV is widely used in image filtering and restoration. However, in this work, to simplify the presentation of the TV filtering algorithm, we focus solely on one-dimensional signal filtering.

In addition to the aforementioned applications, the TV algorithm has a broad capacity for extension to solve various other problems. This includes fields such as image processing, computer vision, signal processing (as discussed in this work), tomography, and medical imaging, as well as being applicable to video processing. This versatility highlights the relevance and utility of this algorithm in different research contexts and practical applications. Recently, faster algorithms for TV filtering have been developed. The development of fast and robust algorithms for nonlinear TV-related filtering is an active area of research.

In the Rudin-Osher-Fatemi image noise removal model \cite{rudin1992nonlinear}, total variation (TV) is employed as a global regularization term. However, we note that the local interactions generated by total variation have a limited range at long distances in practice, placing the model in a position close to a local filter. This observation emphasizes the importance of exploring approaches that can capture local interactions more efficiently, especially in contexts where preserving fine details plays a crucial role, such as in high-resolution image processing applications \cite{louchet2011total}.

Any algorithm that solves the optimization problem can be used to implement TV noise removal. However, the implementation is not trivial because the cost function is not differentiable. Various algorithms have been developed to solve the TV filtering problem, such as convolutional neural network filtering \cite{lan2021image}, bilateral filtering \cite{tomasi1998bilateral}, and non-local filtering \cite{buades2011non}. Total variation is used not only for noise removal but also for more general signal restoration problems, including deconvolution, interpolation, data compression, etc. Moreover, the concept has been generalized and extended in various ways. ///
The central objective of this paper is to develop and optimize filtering algorithms based on the Total Variation technique. These algorithms are specifically applied to one-dimensional signal filtering.

\section{Fundamentals of One-Dimensional Filtering}

\subsection{Signals }
 Signals encountered in electrical engineering and related fields are typically represented as real or complex functions of one or more variables, such as time and/or a variable indicating the outcome of a random experiment \cite{oppenheim2010class}. These signals can be classified in various ways:

\begin{itemize}
    \item One-dimensional or multidimensional
    \item Continuous in time (CT) or discrete in time (DT)
    \item Deterministic or stochastic (also known as random)
\end{itemize}

This diversity of characteristics in signals provides a broad scope for analysis and processing, allowing the application of a variety of techniques and methods to handle different types of signals across a wide range of practical applications.

In this research, I will work with one-dimensional signals that are continuous in time.

A one-dimensional signal is one that varies in only one dimension, whether it be temporal, spatial, or another. This simplified characterization allows for a more direct and specific analysis of the signal's behavior along a single dimension, making it fundamental in a wide range of applications in engineering, data science, and other disciplines \cite{proakis2007digital}.

The understanding of one-dimensional signals is essential for the development of signal processing techniques, mathematical modeling, and statistical analysis in various fields of study and professional practice.

\begin{figure}[H]
    \centering
    \begin{subfigure}{0.5\textwidth}
        \centering
        \includegraphics[width=\linewidth]{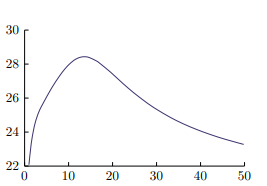}
        \caption{Analog signal in continuous time \cite{getreuer2012rudin}}
        \label{fig:sinal_analogico}
    \end{subfigure}\hfill
    \begin{subfigure}{0.5\textwidth}
        \centering
        \includegraphics[width=\linewidth]{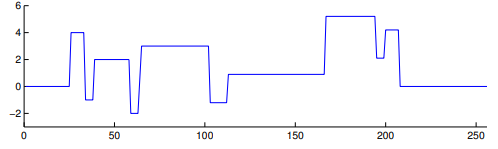}
        \caption{Digital signal in continuous time \cite{selesnick2010total}}
        \label{fig:sinal_digital}
    \end{subfigure}
    \caption{Comparison between signals}
    \label{fig:comparacao_sinais}
\end{figure}

Figures \ref{fig:sinal_analogico} and \ref{fig:sinal_digital} represent analog and digital signals.

\subsection{Filtering}

We assume that we are observing a signal \( x \) that has been corrupted by additive white Gaussian noise, represented by \( y = x + n \), where \( y \), \( x \), and \( n \) belong to the set of real numbers (\( \mathbb{R} \)). In many scenarios, such as wireless communications and signal processing, it is common for signals to be subject to this type of corruption.

To recover the original signal \( x \) from \( y \), a common approach is to use the total variation (TV) filtering technique. This technique is useful for suppressing noise while preserving the edges and important features of the original signal \cite{rudin1992nonlinear}. The objective is to find the signal \( x \) that minimizes the cost function \ref{Função Custo}:

\begin{equation}
    \text{J[x]} = ||y - x||^2_2 + \lambda ||Dx||^1_1 
    \label{Função Custo}
\end{equation}

Where \( ||y - x||^2_2 \) is the \( L^2 \) norm (Euclidean norm) of the difference between \( y \) and \( x \), \( ||Dx||^1_1 \) is the \( L^1 \) norm of the derivative of the signal \( x \), and \( \lambda \) is the regularization parameter that controls the balance between fidelity to the observed signal \( y \) and the smoothing of the resulting signal \( x \). The \( ||Dx||^1_1 \) part of the cost function promotes the smoothing of the signal \( x \), while \( ||y - x||^2_2 \) aims to minimize the error between \( x \) and \( y \).

It is important to note that the regularization parameter \( \lambda \) plays a crucial role in the effectiveness of the TV filtering process. For higher levels of noise in \( y \), a larger value of \( \lambda \) is required to ensure adequate noise suppression and accurate recovery of the original signal \( x \). Therefore, adjusting \( \lambda \) according to the intensity of the noise is essential for obtaining satisfactory results in the TV filtering process.

\subsection{Total Variation Filtering}
Total variation filtering assumes that the noisy data \( y(n) \) takes the form

\[
y(n) = x(n) + w(n), \quad n = 0, \ldots, N - 1
\]

The equation \( y(n) = x(n) + w(n) \) describes the relationship between an observed signal \( y(n) \), a true signal \( x(n) \), and the noise \( w(n) \) in a signal processing system. At a given time \( n \), the observed signal is the combination of the true signal and the noise present at that moment. This relationship is extremely important in various signal processing applications, such as recovering true signals from noisy observations, detecting events amid interference, or analyzing experimental data. Understanding and modeling this relationship is essential for developing effective filtering, reconstruction, or analysis techniques for signals in noisy environments.

We assume that we are observing a signal \( x \) that has been corrupted by additive white Gaussian noise, represented by \( y = x + n \), where \( y \), \( x \), and \( n \) belong to the set of real numbers (\( \mathbb{R} \)). In many scenarios, such as wireless communications and signal processing, it is common for signals to be subject to this type of corruption.

To recover the original signal \( x \) from \( y \), a common approach is to use the total variation (TV) filtering technique. This technique is useful for suppressing noise while preserving the edges and important features of the original signal \cite{rudin1992nonlinear}. The objective is to find the signal \( x \) that minimizes the cost function \ref{Função Custo}:

\[
\text{J[x]} = ||y - x||^2_2 + \lambda ||Dx||^1_1
\]

Where \( ||y - x||^2_2 \) is the \( L^2 \) norm (Euclidean norm) of the difference between \( y \) and \( x \), \( ||Dx||^1_1 \) is the \( L^1 \) norm of the derivative of the signal \( x \), and \( \lambda \) is the regularization parameter that controls the balance between fidelity to the observed signal \( y \) and the smoothing of the resulting signal \( x \). The \( ||Dx||^1_1 \) part of the cost function promotes the smoothing of the signal \( x \), while \( ||y - x||^2_2 \) aims to minimize the error between \( x \) and \( y \).

It is important to note that the regularization parameter \( \lambda \) plays a crucial role in the effectiveness of the TV filtering process. For higher levels of noise in \( y \), a larger value of \( \lambda \) is required to ensure adequate noise suppression and accurate recovery of the original signal \( x \). Therefore, adjusting \( \lambda \) according to the intensity of the noise is essential for obtaining satisfactory results in the TV filtering process.

\subsection{Total Variation Filtering}
Total variation filtering assumes that the noisy data \( y(n) \) takes the form

\[
y(n) = x(n) + w(n), \quad n = 0, \ldots, N - 1
\]

The equation \( y(n) = x(n) + w(n) \) describes the relationship between an observed signal \( y(n) \), a true signal \( x(n) \), and the noise \( w(n) \) in a signal processing system. At a given time \( n \), the observed signal is the combination of the true signal and the noise present at that moment. This relationship is extremely important in various signal processing applications, such as recovering true signals from noisy observations, detecting events amid interference, or analyzing experimental data. Understanding and modeling this relationship is essential for developing effective filtering, reconstruction, or analysis techniques for signals in noisy environments.

The regularization parameter \( \lambda > 0 \) controls the degree of smoothing. Increasing \( \lambda \) gives more weight to the second term, which measures the fluctuation of the signal \( x(n) \) \cite{selesnick2012total}.

The total variation of \( x \) can also be expressed as:
\begin{equation}
\text{TV}(x) = ||Dx||^1_1
\end{equation}
Where \( || * ||_1 \) is the \( L^1 \) norm and

\[
D = \begin{bmatrix}
-1 & 1 & 0 & \cdots & 0 \\
0 & -1 & 1 & \cdots & 0 \\
\vdots & \vdots & \vdots & \ddots & \vdots \\
0 & 0 & \cdots & -1 & 1 \\
\end{bmatrix}
\]

is a matrix of dimensions \( (N - 1) \times N \). Minimizing the equation \ref{Função Custo} is mathematically equivalent to minimizing 
\begin{equation}
J(x) = ||y - x||^2_2 + \lambda ||Dx||^1_1
\label{FunçaoCustoSimplificad}
\end{equation}

This approach is called total variation filtering. Higher noise levels require larger values of \( \lambda \).

The function \( J(x) \) represents the optimization criterion used in the total variation-based filtering technique. This function consists of two main terms. The first term, \( \|y-x\|^2_2 \), measures the squared difference between the observed signal \( y \) and the estimate \( x \), penalizing significant divergences between the two signals. The second term, \( \lambda \|Dx\|^1_1 \), is known as the regularization term and promotes the smoothing of the solution \( x \), where \( \lambda \) is a tuning parameter and \( \|Dx\|^1_1 \) is the \( L_1 \) norm of the first difference of \( x \). This regularization term is crucial to avoid overfitting the observed data and to promote smoother and more coherent solutions, making the approach robust to noise while preserving important features of the original signal. Proper adjustment of the parameter \( \lambda \) is fundamental to balancing fidelity to the data and smoothing of the final solution. Thus, the function \( J(x) \) promotes a \textit{trade-off} between the accuracy of reconstruction and the regularity of the solution, being essential in the formulation and resolution of signal filtering problems using total variation techniques.

\section{Iterative Clipping Algorithm}

The Iterative Shrinkage Algorithm for Total Variation Filtering is a technique used to remove noise from a signal, such as a noisy music signal \cite{selesnick2012total}. It starts with the noisy signal and iteratively adjusts a vector to smooth the signal. Each iteration calculates how much the signal differs from the original and applies a filter to reduce these differences, keeping the signal as close as possible to the original.

To do this, the algorithm uses two main ingredients: the error between the original and the adjusted signal, and a penalty for abrupt changes in the signal. It gradually adjusts the filter to balance noise removal without losing important information from the signal. The goal is to find a balance that makes the sound clearer by removing noise without distorting the original signal too much.

Based on what has been said, the signal filtering process is described in Pseudocode \ref{algoritmoRecorte}:

\begin{algorithm}
\caption{Iterative Shrinkage Algorithm for Total Variation Filtering}
\KwIn{Noisy signal $y$, regularization parameter $\lambda$, number of iterations Nit}
\KwOut{Denoised signal $x$, objective function $J$}
Initialize $J$ as a zero vector with length Nit\;
Initialize $N$ as the length of $y$\;
Initialize $z$ as a zero vector with length $N-1$\;
Set $\alpha = 4$\;
Set $T = \frac{\lambda}{2}$\;
\For{$k = 1$ \KwTo Nit}{
    Compute $x = y - [-z(1), -\text{diff}(z), z(\text{end})]$ \;
    Compute $J(k) = \sum |x-y|^2 + \lambda \sum |x_{i+1} - x_i|$\;
    Update $z = z + \frac{1}{\alpha} \, \text{diff}(x)$ 
    Clip each element of $z$ to be within $[-T, T]$\;
}
\label{algoritmoRecorte}
\end{algorithm}

In this context:
The function \texttt{diff} in many programming languages is used to calculate the difference between consecutive elements of a vector or array.
$T$ is defined as half of the regularization parameter $\lambda$. $\alpha$ is a scaling factor.

The algorithm iterates over the estimate of the clean signal $x$ and adjusts the vector $z$ to minimize the objective function, which combines the error between $x$ and $y$ with a penalty for significant variations between adjacent elements in $x$. The clipping in $z$ limits the magnitude of the changes allowed in each iteration.

\section{Results and Discussions}
\label{sec:resultados}

In this section, we discuss the results obtained after applying the filtering algorithm based on the Lagrange multiplier method, emphasizing the results obtained through the analysis of the L curve.

\subsection{Step Signal with Gaussian Noise}
To evaluate the effectiveness of the proposed method, we added 10\% Gaussian noise to a step signal and applied the noise removal algorithm.

\begin{figure}[H]
    \centering
    \begin{subfigure}[b]{0.45\linewidth}
        \centering
        \includegraphics[width=\linewidth]{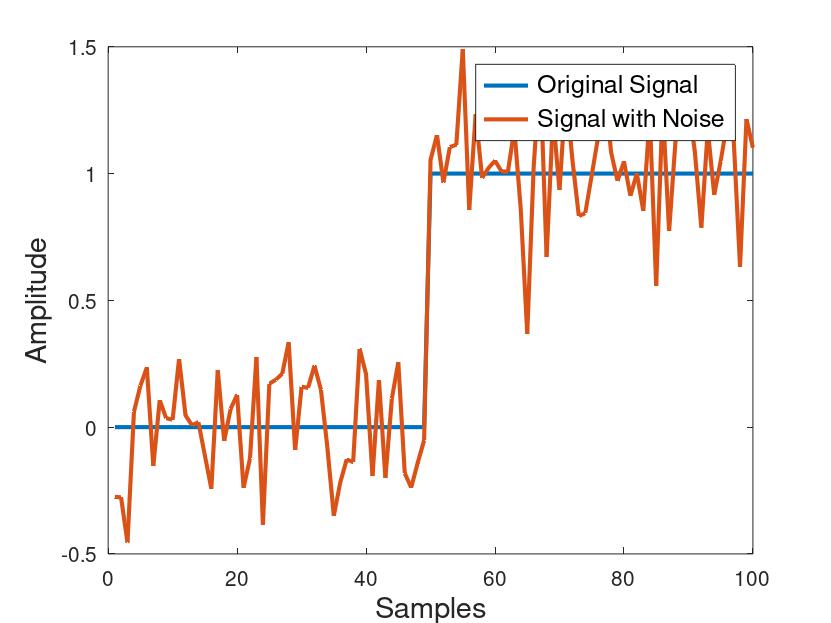}
        \caption{Step Signal with Gaussian Noise}
        \label{fig:sinalDegrau}
    \end{subfigure}
    \hfill
    \begin{subfigure}[b]{0.45\linewidth}
        \centering
        \includegraphics[width=\linewidth]{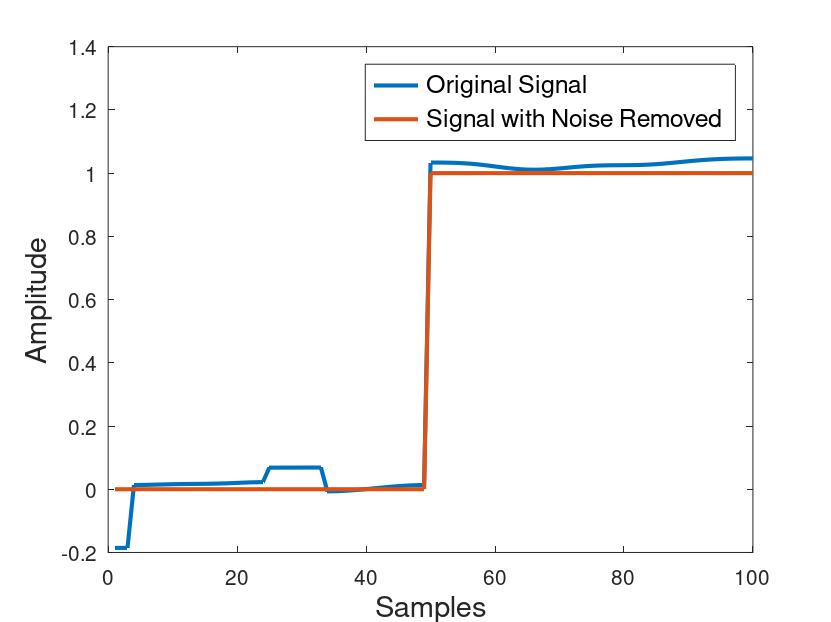}
        \caption{Step Signal after Filtering}
        \label{fig:sinalDegrauRemovido}
    \end{subfigure}
    \caption{Processing step signals: (a) Step Signal with Gaussian Noise and (b) Step Signal after Filtering.}
    \label{fig:degrauProcessamento}
\end{figure}

Initially, we observe the step signal with added Gaussian noise, as shown in Figure \ref{fig:sinalDegrau}.

\begin{figure}[H]
    \centering
    \includegraphics[width=0.5\linewidth]{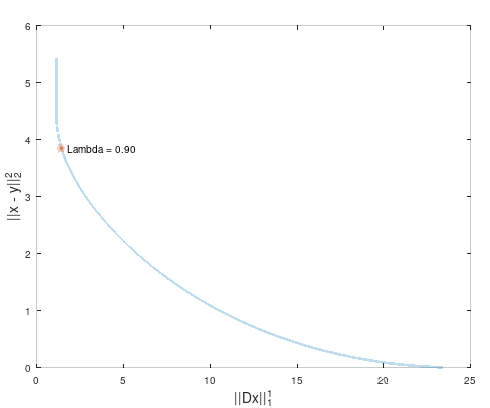}
    \caption{L Curve for the Step Signal}
    \label{fig:curvaLDegrau}
\end{figure}

After applying the filtering algorithm with the Lagrange multiplier adjusted to $\alpha = 0.9$, as determined by the analysis of the L curve (Figure \ref{fig:curvaLDegrau}), we obtained the restored step signal, as shown in Figure \ref{fig:sinalDegrauRemovido}.

The use of the L curve allowed us to determine the parameter $\alpha = 0.9$, significantly optimizing the noise removal while preserving the fundamental characteristics of the original signal. The visual comparison between Figures \ref{fig:sinalDegrau} and \ref{fig:sinalDegrauRemovido} clearly demonstrates the effectiveness of the proposed method in restoring the signal degraded by Gaussian noise.

\subsection{Laplace Signal with Gaussian Noise}

In the second example, we applied the filtering algorithm to a square wave signal with Laplace noise, exploring its applicability to more complex signals.

\begin{figure}[H]
    \centering
    \begin{subfigure}[b]{0.7\linewidth}
        \centering
        \includegraphics[width=\linewidth]{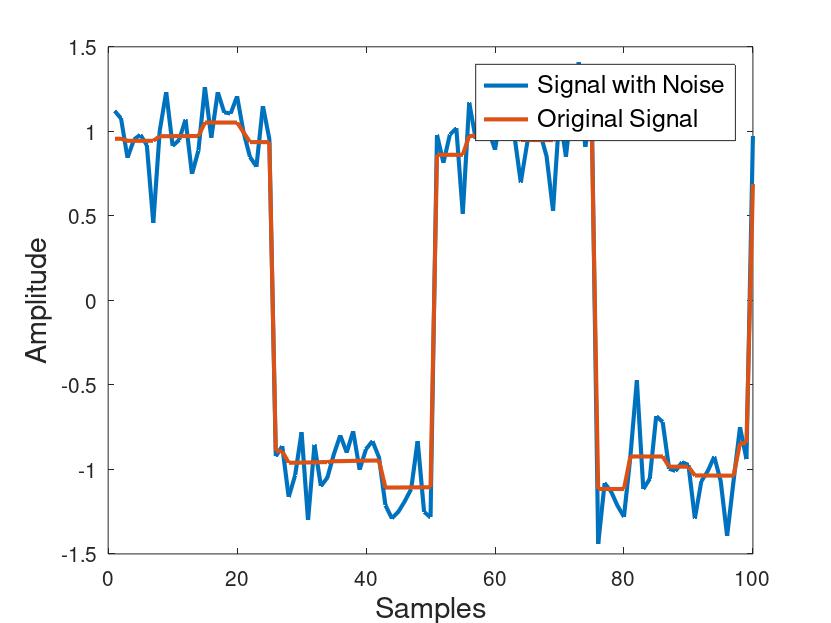}
        \caption{Laplace Signal with Gaussian Noise}
        \label{fig:sinalLaplace}
    \end{subfigure}
    \hfill
    \begin{subfigure}[b]{0.7\linewidth}
        \centering
        \includegraphics[width=\linewidth]{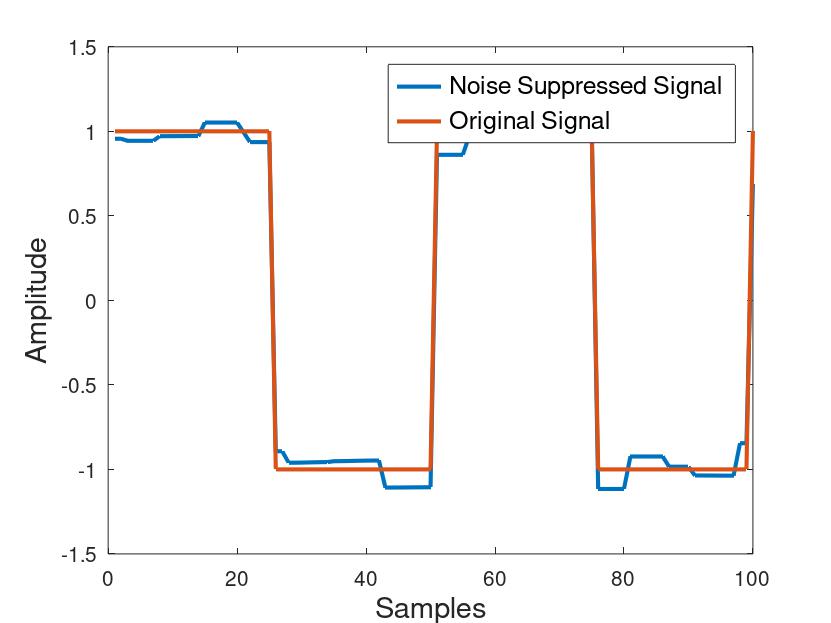}
        \caption{Laplace Signal after Filtering}
        \label{fig:sinalLaplacePosAlgoritmo}
    \end{subfigure}
    \caption{Processing Laplace signals: (a) Laplace Signal with Gaussian Noise and (b) Laplace Signal after Filtering.}    \label{fig:laplaceProcessamento}
\end{figure}

Figure \ref{fig:sinalLaplace} presents the Laplace signal with added Gaussian noise.

\begin{figure}[H]
    \centering
    \caption{L Curve for the Laplace Signal}
    \includegraphics[width=0.8\linewidth]{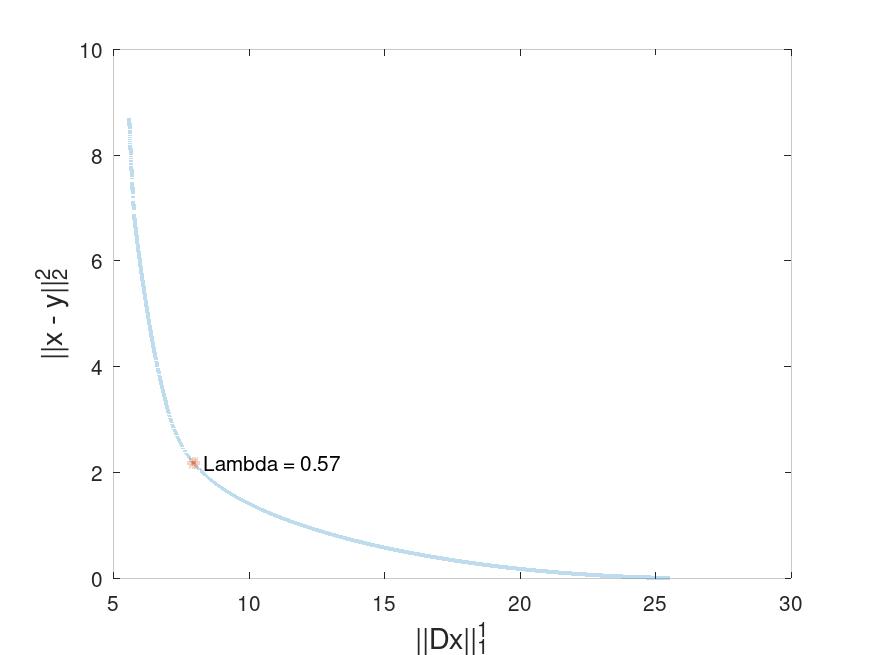}
    \label{fig:curvaLLaplace}
\end{figure}

After applying the noise removal algorithm with the parameter $\alpha = 0.57$, determined by the analysis of the L curve (Figure \ref{fig:curvaLLaplace}), we obtained the restoration of the Laplace signal, as seen in Figure \ref{fig:sinalLaplacePosAlgoritmo}.
The detailed analysis of the L curve provided a precise adjustment of the parameter $\alpha = 0.57$, resulting in a remarkable reduction in noise while preserving the essential characteristics of the original Laplace signal. Figures \ref{fig:sinalLaplace} and \ref{fig:sinalLaplacePosAlgoritmo} demonstrate the effectiveness of the algorithm in improving the signal quality after filtering.

\subsection{Noise Removal from the Vuvuzela}
Finally, we explored the application of the algorithm for removing the characteristic noise of the vuvuzela, a common challenge in sporting events \cite{itu2010}. This study is particularly relevant due to the negative impact that the continuous and deafening sound of vuvuzelas can have on the auditory experience of spectators and the clarity of narration during games.

\begin{figure}[H]
    \centering
    \begin{subfigure}[b]{0.45\linewidth}
        \centering
        \includegraphics[width=\linewidth]{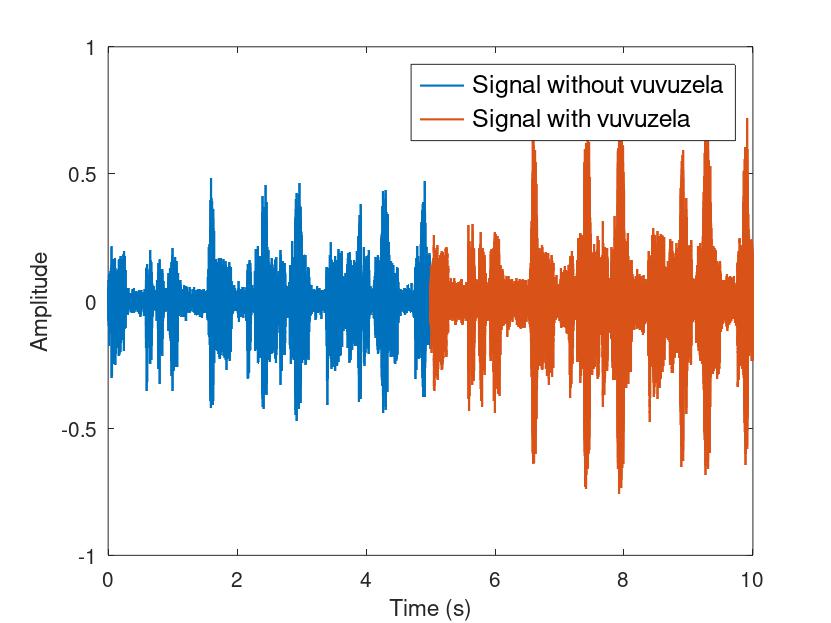}
        \caption{Vuvuzela Sound without Noise}
        \label{fig:sinalVuvuzela}
    \end{subfigure}
    \hfill
    \begin{subfigure}[b]{0.45\linewidth}
        \centering
        \includegraphics[width=\linewidth]{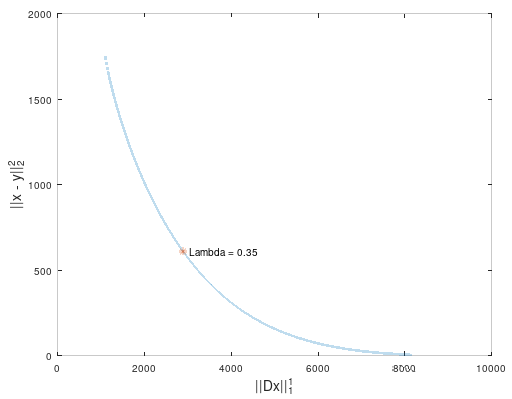}
        \caption{L Curve for the Vuvuzela Sound}
        \label{fig:CurvaLVuvuzela}
    \end{subfigure}
    \caption{Results of vuvuzela noise removal: (a) Vuvuzela Sound without Noise and (b) L Curve for the Vuvuzela Sound.}
    \label{fig:vuvuzelaResultados}
\end{figure}

Figure \ref{fig:sinalVuvuzela} presents a 30-second excerpt of audio from a Brazilian national team game, extracted from the YouTube platform, where the sound of the vuvuzela with noise is clearly audible. This excerpt illustrates the problem we aim to solve, highlighting the need for effective noise removal techniques.

After applying the noise removal algorithm with the parameter $\alpha = 0.35$, adjusted based on the L curve (Figure \ref{fig:CurvaLVuvuzela}), we achieved effective removal of the vuvuzela noise. The results can be observed in Figure \ref{fig:sinalVuvuzela}.

The analysis of the L curve was fundamental in determining the parameter $\alpha = 0.35$, optimizing the removal of vuvuzela noise while preserving the integrity of the audio. Figure \ref{fig:sinalVuvuzela} highlights the effectiveness of the algorithm in significantly reducing the characteristic noise of the vuvuzela.

These results demonstrate the applicability and effectiveness of the algorithm based on the Lagrange multiplier method for noise removal in different signal processing contexts. The use of the L curve as an adjustment tool provided significant improvements in the quality of processed signals, validating its efficiency and relevance for practical applications.

To listen to the original audio with the vuvuzela sound, visit \href{https://github.com/Joyce8900/TCC}{link}. To listen to the results obtained from the removal of the vuvuzela noise, visit \href{https://github.com/Joyce8900/TCC}{link}.

\section{Conclusion}
In summary, this work demonstrates that the proposed algorithm excels in solving complex noise removal problems. However, it is essential to recognize that, in some of the studied cases, the results did not achieve 100\% accuracy. This variation in performance can be attributed to the complexity of the data, the limitations of the adopted model, and the specific conditions of each dataset.

In the specific example of vuvuzela noise removal, although we observed significant improvements in audio quality, the presence of the noise frequency often interfered with the frequency of the sports commentator's voice, preventing complete noise elimination. This interference highlights additional challenges in scenarios where the noise and the signal of interest share similar spectral characteristics.

Despite these variations, the results obtained reinforce the robustness and effectiveness of the algorithm in manipulating noisy signals, emphasizing its potential for various practical applications. The continuous evolution and refinement of the method, coupled with personalized adjustments based on the analysis of the L curve, promise significant advancements in addressing complex signal processing challenges in the future.

\bibliographystyle{unsrt}  
\bibliography{references}  

\end{document}